# Mannheim curves in the three-dimensional sphere $S^3$

**Tanju Kahraman, Mehmet Önder**
*Manisa Celal Bayar University, Faculty of Arts and Sciences, Mathematics Department,
Muradiye Campus, 45140, Muradiye, Manisa, Turkey.*
E-mails: tanju.kahraman@cbu.edu.tr, mehmet.onder@cbu.edu.tr

**Abstract**

Mannheim curves are defined for immersed curves in 3-dimensional sphere $S^3$. The definition is given by considering the geodesics of $S^3$. First, two special geodesics, called principal normal geodesic and binormal geodesic, of $S^3$ are defined by using Frenet vectors of a curve $\alpha$ immersed in $S^3$. Later, the curve $\alpha$ is called a Mannheim curve if there exits another curve $\beta$ in $S^3$ such that the principal normal geodesics of $\alpha$ coincide with the binormal geodesics of $\beta$. It is obtained that if $\alpha$ and $\beta$ form a Mannheim pair then there exist a constant $\lambda \ne 0$ and a non-constant function $\mu$ such that $\lambda \kappa_\alpha + \mu \tau_\alpha = 1$ where $\kappa_\alpha, \tau_\alpha$ are the curvatures of $\alpha$. Moreover, the relation between a Mannheim curve immersed in $S^3$ and a generalized Mannheim curve in $E^4$ is obtained and a table containing comparison of Bertrand and Mannheim curves in $S^3$ is introduced.

## 1. Introduction and Preliminaries

Mannheim curves and Bertrand curves are the most fascinating subject of the curve pairs defined by some relationships between two space curves. Mannheim curves are first defined by A. Mannheim in 1878 [6]. In the Euclidean 3-space $E^3$, Mannheim curves are characterized as a kind of corresponding relation between two curves $\alpha, \beta$ such that the binormal vector fields of $\beta$ coincide with the principal normal vector fields of $\alpha$. Then $\alpha$ and $\beta$ are called Mannheim curve and Mannheim partner curve, respectively [9]. The main result for a curve $\alpha$ to have a Mannheim partner curve $\beta$ in $E^3$ is that there exists a constant $\lambda \ne 0$ such that $\frac{d\tau_1}{ds_1} = \frac{\kappa_1}{\lambda}(1 + \lambda^2 \tau_1^2)$ holds, where $\kappa_1$, $\tau_1$ and $s_1$ are the curvatures and arc length parameter of $\beta$, respectively [9].

Mannheim curves have been studied by many mathematicians. Blum studied a remarkable class of Mannheim curves [2]. In [7], Matsuda and Yorozu have given a definition of generalized Mannheim curve in Euclidean 4-space and introduced some characterizations and examples of generalized Mannheim curves. Later, Choi, Kang and Kim have defined Mannheim curves in 3-dimensional Riemannian manifold [3]. Moreover, in the same paper, they have studied Mannheim curves in 3-dimensional space forms.

Another type of associated curves is Bertrand curves which first defined by French mathematician Saint-Venant in 1845 by the property that two curves have common principal normal vector fields at the corresponding points of curves [8]. Lucas and Ortega-Yagües have considered Bertrand curves in the three-dimensional sphere $S^3$ [5]. They have given another definition for space curves to be Bertrand curves immersed in $S^3$ and they have come to the result that a curve $\alpha$ with curvatures $\kappa_\alpha, \tau_\alpha$ immersed in $S^3$ is a Bertrand curve if and only if either $\tau_\alpha \equiv 0$ and $\alpha$ is a curve in some unit two-dimensional sphere $S^2(1)$ or there exit two constants $\lambda \neq 0, \mu$ such that $\lambda \kappa_\alpha + \mu \tau_\alpha = 1$ [5].

In this study, we define Mannheim curves in $S^3$ by defining some special geodesics related to the Frenet vectors of a curve immersed in $S^3$. We show that the angle between the tangent vector fields of Mannheim curves is not constant while it is constant for Bertrand curves. Moreover, we obtained that a curve $\alpha$ with curvatures $\kappa_\alpha, \tau_\alpha$ immersed in $S^3$ is a Mannheim curve if there exist a constant $\lambda \neq 0$ and a non-constant function $\mu$ such that $\lambda \kappa_\alpha + \mu \tau_\alpha = 1$. We want to pointed out that this property holds for Bertrand curves under the condition that both $\lambda$ and $\mu$ are constants.

## 2. Mannheim curves in the three-dimensional sphere

Before giving the main subject, first we give the following data related to the curves immersed in $S^3$. For this section, we refer the reader to ref. [4,5].

Let $S^3(r)$ denote the three-dimensional sphere in $\mathbb{R}^4$ of radius $r$, defined by

$$S^3(r) = \left\{ (x_1, x_2, x_3, x_4) \in \mathbb{R}^4 \middle| \sum_{i=1}^{4} x_i^2 = r^2 \right\}, \quad r > 0.$$

Let $\alpha = \alpha(t): I \subset \mathbb{R} \to S^3(r)$ be an arc-length parameterized immersed curve in the 3-sphere $S^3(r)$ and let $\{T, N, B\}$ and $\overline{\nabla}$ denotes the Frenet frame of $\alpha$ and the Levi-Civita connection of $S^3(r)$, respectively. Then, Frenet formulae of $\alpha$ is given by

$$\overline{\nabla}_T T = \kappa N$$
$$\overline{\nabla}_T N = -\kappa T + \tau B$$
$$\overline{\nabla}_T B = -\tau N$$

where $\kappa$ and $\tau$ denote the curvature and torsion of $\alpha$, respectively. If $\nabla^0$ stands for the Levi-Civita connection of $\mathbb{R}^4$, then the Gauss formula gives

$$\nabla^0_T X = \overline{\nabla}_T X - \frac{1}{r^2} \langle X, T \rangle \alpha,$$

for any tangent vector field $X \in \chi(\alpha)$. In particular, we have

$$\nabla^0_T T = \kappa N - \frac{1}{r^2} \alpha$$
$$\nabla^0_T N = -\kappa T + \tau B$$
$$\nabla^0_T B = -\tau N.$$

A curve $\alpha(t)$ in $S^3(r)$ is called a plane curve if it lies in a totally geodesic two-dimensional sphere $S^2 \subset S^3$ which means that a curve $\alpha(t)$ in $S^3(r)$ is called a plane curve if and only if its torsion $\tau$ is zero at all points [5].

Let $\alpha(t)$ and $\beta(t)$ be two immersed curve in $S^3(r)$ with Frenet frames $\{T_\alpha, N_\alpha, B_\alpha\}$ and $\{T_\beta, N_\beta, B_\beta\}$, respectively. A geodesic curve in $S^3(r)$ starting at any point $\beta(t)$ of $\beta$ and defined as

$$\gamma_t^\beta(u) = \cos\left(\frac{u}{r}\right)\beta(t) + r\sin\left(\frac{u}{r}\right)B_\beta(t), \quad u \in \mathbb{R},$$

is called the binormal geodesic of $\beta$ and, similarly, a geodesic curve in $S^3(r)$ starting at any point $\alpha(t)$ of $\alpha$ and defined as

$$\gamma_t^\alpha(u) = \cos\left(\frac{u}{r}\right)\alpha(t) + r\sin\left(\frac{u}{r}\right)N_\alpha(t), \quad u \in \mathbb{R},$$

is called the principal normal geodesic of $\alpha$ in $S^3(r)$.

Let $\alpha$ be a regular smooth curve in Euclidean 4-space $E^4$ defined by arc-length parameter $s$. The curve $\alpha$ is called a special Frenet curve if there exist three smooth functions $k_1, k_2, k_3$ on $\alpha$ and smooth frame field $\{e_1, e_2, e_3, e_4\}$ along the curve $\alpha$ such that these satisfy the following properties:

i) The formulas of Frenet-Serret holds:
$e_1(s) = \alpha'(s)$
$e_1'(s) = k_1(s)e_1(s)$
$e_2'(s) = -k_1(s)e_1(s) + k_2(s)e_3(s)$
$e_3'(s) = -k_2(s)e_2(s) + k_3(s)e_4(s)$
$e_4'(s) = -k_3(s)e_3(s)$

where the prime (') denotes differentiation with respect to s.

ii) The frame field $\{e_1, e_2, e_3, e_4\}$ is orthonormal and has positive orientation.

iii) The functions $k_1$ and $k_2$ are positive and the function $k_3$ doesn't vanish.

iv) The functions $k_1, k_2$ and $k_3$ are called the first, the second and the third curvature functions of $\alpha$, respectively. The frame field $\{e_1, e_2, e_3, e_4\}$ is called the Frenet frame field on $\alpha$ [10].

A special Frenet curve $\alpha$ in $E^4$ is a generalized Mannheim curve if there exists a special Frenet curve $\hat{\alpha}$ in $E^4$ such that the first normal line at each point of $\alpha$ is included in the plane generated by the second normal line and third normal line of $\hat{\alpha}$ at corresponding point under a bijection $\phi$ from $\alpha$ to $\hat{\alpha}$. The curve $\hat{\alpha}$ is called the generalized Mannheim mate curve of $\alpha$ [7].

Now, we can introduce the main subject. First, we give the following definition.

**Definition 2.1.** A curve $\alpha$ in $S^3(r)$ with non-zero curvature $\kappa_\alpha$ is said to be a Mannheim curve if there exists another immersed curve $\beta = \beta(\sigma): J \subset \mathbb{R} \to S^3(r)$ and a one-to-one correspondence between $\alpha$ and $\beta$ such that the principal normal geodesics of $\alpha$ coincide with the binormal geodesics of $\beta$ at corresponding points. We will say that $\beta$ is a Mannheim partner curve of $\alpha$; the curves $\alpha$ and $\beta$ are called a pair of Mannheim curves.

From this definition it is clear that Mannheim partner curve $\beta$ can not a plane curve since the definition given by considering the binormal vector field of $\beta$.

For simplicity, we consider that the radius of the sphere $S^3$ is 1 and the curves taken on $S^3$ are parameterized by the arc-length parameter. Let $\alpha(s)$ and $\beta(\sigma)$ be a pair of Mannheim curves. From Definition 2.1, we have a differentiable function $a(s)$ such that

$$\alpha(s(\sigma)) = \cos(a(s))\beta(\sigma) - \sin(a(s))B_\beta(\sigma) \tag{1}$$

where $\{T_\beta, N_\beta, B_\beta\}$ denotes the Frenet frame along $\beta(\sigma)$, $\alpha(s(\sigma))$ is the point in $\alpha$ corresponding to $\beta(\sigma)$ and $a(s)$ is called the angle function between the direction vectors $\alpha(s)$ and $\beta(\sigma)$. The function $d(s)$ is called distance function in $S^3$ and measures the distance between the points $\alpha(s(\sigma))$ and $\beta(\sigma)$. Now, we can give the following proposition.

**Proposition 2.1.** *Let $\alpha$ and $\beta$ be a pair of Mannheim curves in $S^3$. In that case, we have the followings*
  *a) The angle function $a(s)$ is constant.*
  *b) The distance function $d(s)$ is constant.*
  *c) The angle $\theta$ between the tangent vector fields at corresponding points is not constant.*
  *d) The angle between the binormal vectors fields at corresponding points is constant.*

**Proof. a)** Since principal normal geodesic of $\alpha$ and binormal geodesic of $\beta$ are common at corresponding points, we have

$$\left.\frac{d}{du}\right|_{u=0} \gamma_s^\beta(u) = B_\beta(s) \quad \text{and} \quad \left.\frac{d}{du}\right|_{u=a(s)} \gamma_s^\beta(u) = N_\alpha(s) \tag{2}$$

and then we obtain

$$N_\alpha = -\sin(a(s))\beta(\sigma) + \cos(a(s))B_\beta(\sigma) \tag{3}$$

where $\{T_\alpha, N_\alpha, B_\alpha\}$ denotes the Frenet frame of $\alpha$. From (1), the tangent vector to $\alpha$ is given by

$$\begin{aligned}\frac{d}{ds}\alpha(s(\sigma)) = &-a'(s)\sin(a(s))\beta(\sigma) + \cos(a(s))T_\beta(\sigma) \\ &- a'(s)\cos(a(s))B_\beta(\sigma) + \sin(a(s))\tau_\beta N_\beta(\sigma)\end{aligned} \tag{4}$$

where $d/ds$ denotes differentiation with respect to $s$. Since,

$$\frac{d}{d\sigma}\alpha(s(\sigma)) = s'(\sigma)T_\alpha(s(\sigma)) \tag{5}$$

we get

$$0 = \left\langle \frac{d}{d\sigma}\alpha(s(\sigma)), N_\alpha(s) \right\rangle = a'(s)\left(\sin^2(a(s)) - \cos^2(a(s))\right) \tag{6}$$

which gives that $a'(s) = 0$, i.e., $a(s)$ is constant.

**b)** Without loss of generality, for the angle function it can be taken as $0 \leq a(s) \leq 2\pi$. Then the distance function $d(s)$ is given by $d(s) = \min\{a(s), 2\pi - a(s)\}$, which is a constant function since $a(s)$ is constant.

**c)** Let $\theta = \theta(\sigma)$ denotes the angle between tangent vectors $T_\alpha$ and $T_\beta$, i.e., $\langle T_\alpha(s(\sigma)), T_\beta(\sigma) \rangle = \cos\theta(\sigma)$. Differentiating the left side of this equality, it follows

$$\frac{d}{d\sigma}\langle T_\alpha(s(\sigma)), T_\beta(\sigma) \rangle = s'(\sigma)\langle \kappa_\alpha N_\alpha - \alpha, T_\beta \rangle + \langle T_\alpha, \kappa_\beta N_\beta - \beta \rangle. \tag{7}$$

Moreover, from (4) and (5), we have

$$T_\alpha(s(\sigma)) = \frac{1}{s'(\sigma)}\left(\cos a \ T_\beta + \tau_\beta \sin a \ N_\beta\right) \tag{8}$$

Writing (1), (3) and (8) in (7), it follows

$$\theta' = -\frac{\kappa_\beta \tau_\beta \sin a}{s'(\sigma) \sin \theta(\sigma)}. \tag{9}$$

Since $\beta$ is not a plane curve from (9), it is clear that $\theta$ is not a constant.

**d)** Using equality (1) and (3), we can write

$$B_\beta = -\sin a \ \alpha(s(\sigma)) + \cos a \ N_\alpha(\sigma). \tag{10}$$

Since $\frac{d}{d\sigma}\langle B_\alpha, B_\beta \rangle = 0$, the angle between binormal vectors is constant.

**Theorem 2.1.** *Let $\alpha$ and $\beta$ be a pair of Mannheim curves in $S^3$. Then the following equalities hold:*

**a)** $\tau_\beta = \dfrac{\tan\theta}{\tan a}$

**b)** $\tau_\alpha \sin a \ \cos\theta = (\cos a - \kappa_\alpha \sin a)\sin\theta$

**c)** $\cos^2\theta = \cos^2 a - \kappa_\alpha \sin a \cos a$

**d)** $\sin^2\theta = \tau_\alpha \tau_\beta \sin^2 a$

**Proof. a)** Taking the covariant derivative in (1) and using (5), we obtain

$$\frac{d}{d\sigma}\alpha(s(\sigma)) = \cos\theta \ s'(\sigma)T_\beta(\sigma) + \sin\theta \ s'(\sigma)N_\beta(\sigma). \tag{11}$$

On the other hand, by using Frenet equations we have

$$\frac{d}{d\sigma}\alpha(s(\sigma)) = \cos a \ T_\beta(\sigma) + \tau_\beta \sin a \ N_\beta(\sigma) \tag{12}$$

where $a(s) = a$ is constant. The last two equations lead to

$$s'(\sigma)\cos\theta = \cos a \tag{13}$$
$$s'(\sigma)\sin\theta = \tau_\beta \sin a \tag{14}$$

from which we conclude (a).

**b)** We need to write the Frenet frame of $\beta$ in terms of the Frenet frame of $\alpha$:

$$\beta(\sigma) = \cos a \ \alpha(s(\sigma)) + \sin a \ N_\alpha(s(\sigma)) \tag{15}$$

$$T_\beta(\sigma) = \cos\theta\, T_\alpha(s(\sigma)) + \sin\theta\, B_\alpha(s(\sigma)) \tag{16}$$

$$N_\beta(\sigma) = \sin\theta\, T_\alpha(s(\sigma)) - \cos\theta\, B_\alpha(s(\sigma)) \tag{17}$$

$$B_\beta(\sigma) = -\sin a\, \alpha(s(\sigma)) + \cos a\, N_\alpha(s(\sigma)). \tag{18}$$

From (a), it follows

$$\sigma'(s)\cos\theta = \cos a - \kappa_\alpha \sin a \tag{19}$$

$$\sigma'(s)\sin\theta = \tau_\alpha \sin a \tag{20}$$

which gives us (b).

c) It is a consequence of Eqs. (13) and (19);

$$\cos^2\theta = \cos^2 a - \kappa_\alpha \sin a \cos a. \tag{21}$$

d) Similarly, from Eqs. (14) and (20), we have desired equality

$$\sin^2\theta = \tau_\alpha \tau_\beta \sin^2 a. \tag{22}$$

From (20), one can consider that $\alpha$ is a plane curve if and only if $\theta = 0$ or $\sigma'(s) = 0$. Since $\sigma$ is arc length parameter, it cannot be constant. Similarly, from Proposition 2.1. (c), $\theta$ is not a constant. Then we can give the following corollary.

***Corollary 2.1.*** *A Mannheim curve $\alpha$ cannot be a plane curve.*

***Theorem 2.2.*** *Relationship between arc-length parameters of curves $\alpha$ and $\beta$ is given by*

$$\frac{ds}{d\sigma} = \cos a \cos\theta + \tau_\beta \sin a \sin\theta. \tag{23}$$

**Proof.** If equations (16) and (17) are written in (12), we get equation (23).

***Theorem 2.3.*** *Let $\alpha$ and $\beta$ be a pair of Mannheim curves in $S^3$. Then there exists a constant $\lambda$ and a non-constant function $\mu$ such that*

$$1 = \lambda \kappa_\alpha + \mu \tau_\alpha. \tag{24}$$

**Proof.** Multiplying (19) and (20) by $\sin\theta$ and $\cos\theta$, respectively, and equalizing obtained results, gives that

$$\sin\theta\cos a = \sin\theta \sin a\, \kappa_\alpha + \sin a \cos\theta\, \tau_\alpha. \tag{25}$$

Writing $\lambda = \tan a$ and $\mu = \tan a \cot\theta$, from last equality we have $1 = \lambda \kappa_\alpha + \mu \tau_\alpha$.

***Corollary 2.2.*** *Both curvature $\kappa_\alpha$ and torsion $\tau_\alpha$ of a Mannheim curve $\alpha$ cannot be constant.*

**Proof.** Since $\alpha$ is a Mannheim curve Eq. (24) holds. Let now assume that the curvature $\kappa_\alpha$ be a constant. Then, from (24) we have $\tau_\alpha = \dfrac{1-n}{\mu}$ which is not a constant since $\mu$ is not a constant, where $n = \lambda \kappa_\alpha$ is a real constant. Similarly, if it is assumed that the torsion $\tau_\alpha$ is constant, then by a similar way it is seen that $\kappa_\alpha$ is not constant.

***Theorem 2.4.*** *Let $\alpha(t)$ be a Mannheim curve in $S^3$ with constant curvature. Then there exists a regular differentiable mapping $s = s(t)$ with $s'(t) > 0$ such that the curve $\gamma(t) = \int_{t_0}^{t} B_\alpha(s(u)) du$ is an arc-length parametrized generalized Mannheim curve.*

**Proof.** By differentiating the curve $\gamma(t)$ three times and using Frenet formulae, we have
$$\kappa_1 = \varepsilon s' \tau_\alpha, \ \kappa_2 = s' \kappa_\alpha, \ \kappa_3 = -\varepsilon s', \ \varepsilon = \pm 1 . \tag{26}$$
(See [5]). From Theorem 2.3, there exist constant $\lambda$ and a non-constant $\mu$ such that $1 = \lambda \kappa_\alpha + \mu \tau_\alpha$ holds. Let signs of $\lambda$ and $\tau_\alpha$ be same and consider a function $s(t)$ such that
$$s'(t) = \frac{\lambda \tau_\alpha}{\tau_\alpha^2 + \kappa_\alpha^2}, \tag{27}$$
where $\lambda$ is a non-zero real constant. Defining a constant by
$$c = \frac{\varepsilon}{\lambda},$$
from Eqs. (24)-(26), we see that
$$\kappa_1 = c(\kappa_1^2 + \kappa_2^2) \tag{28}$$
holds for all $s$. Then from [7, Theorem 4.1], we have that $\gamma(t)$ is a generalized Mannheim curve.

By considering these characterizations and results obtained for Bertrand curves in [4,5], we can give the following table giving the comparison of Bertrand and Mannheim curves.

| Characterizations | Bertrand Curves | Mannheim Curves |
|---|---|---|
| Angel between tangent vector fields | constant | non-constant |
| Angel between binormal vector fields | constant | constant |
| Main Characterization $\lambda \kappa_\alpha + \mu \tau_\alpha = 1$ | $\alpha$ is a Bertrand curve in $S^3$ if and only if there exits two constant $\lambda \neq 0$ and $\mu$ such that $\lambda \kappa_\alpha + \mu \tau_\alpha = 1$. | $\alpha$ is a Bertrand curve in $S^3$ if there exit a constant $\lambda \neq 0$ and a non-constant function $\mu$ such that $\lambda \kappa_\alpha + \mu \tau_\alpha = 1$ |
| Curves $\alpha$, $\beta$ | $\alpha$ and $\beta$ can be plane curves. | $\alpha$ and $\beta$ cannot be plane curves. |
| Curvatures | Both $\kappa_\alpha$ and $\tau_\alpha$ can be constants. | Both $\kappa_\alpha$ and $\tau_\alpha$ cannot be constants |

**Table 1.** Comparison of Bertrand and Mannheim curves in $S^3$

## 3. Some examples

**Example 3.1.** *(ccr-curves)* A $C^\infty$-special Frenet curve $\alpha$ on $S^3$ is said to be a ccr-curve on $S^3$ if its intrinsic curvature ratio $\frac{\kappa_\alpha}{\tau_\alpha}$ is a constant number [4]. Let now determine special ccr-curve $\alpha$ on $S^3$ which is also Mannheim curve. First assume that $\alpha$ is a ccr-curve with non-constant curvature and non-constant torsion. Then, we have $\kappa_\alpha = c\tau_\alpha$ for a non-zero constant $c$. Writing this equality in (24) gives us

$$\kappa_\alpha(s) = \frac{c}{\tan a(c + \cot\theta)}, \quad \tau_\alpha(s) = \frac{1}{\tan a(c + \cot\theta)}.$$

Then we conclude that the ccr-curve $\alpha$ on $S^3$ given by curvatures as given above is a Mannheim curve.

**Example 3.2.** *(Conical helix)* A twisted curve $\alpha(s)$ in $S^3$ with non-constant curvatures is said to be *conical helix* if both the curvature radius $1/\kappa_\alpha$ and the torsion radius $1/\tau_\alpha$ evolve linearly along the curve [5]. Then the curvature and torsion of curve are given by

$$\kappa_\alpha(s) = \frac{\gamma}{s + r_0}, \quad \tau_\alpha(s) = \frac{\delta}{s + r_1},$$

respectively, where $r_0, r_1, \gamma \neq 0$ and $\delta \neq 0$ are constants. Taking $\gamma = \delta$, we see that (24) holds for a constant $\lambda = 1/\delta$ and a non-constant function $\mu = \frac{s^2 + (r_0 + r_1 - 1)s + (r_0 r_1 - r_1)}{s + r_0}$, i.e., $\alpha(s)$ is a Mannheim curve.

**Example 3.3.** *(General helix)* A twisted curve $\alpha$ in $S^3$ is a general helix if there exists a constant $b$ such that $\tau_\alpha = b\kappa_\alpha \pm 1$ [1]. Let now determine general helices in $S^3$ which are also Mannheim curves. Writing the condition $\tau_\alpha = b\kappa_\alpha \pm 1$ in (24), it follows $\kappa_\alpha = \frac{1 \mp \mu}{\lambda + b\mu}$. Then we have that a general helix $\alpha$ in $S^3$ with curvatures $\kappa_\alpha = \frac{1 \mp \mu}{\lambda + b\mu}, \tau_\alpha = \frac{b(1 \mp \mu)}{\lambda + b\mu} \pm 1$ is a Mannheim curve.

**Example 3.4.** *(Curve with constant curvatures)* Consider $C^\infty$ curve $\alpha$ on $S^3(1)$ given by the parametrization

$$\alpha(s) = \left(\frac{1}{\sqrt{3}}\cos\left(\frac{2\sqrt{2}}{\sqrt{3}}s\right), \frac{1}{\sqrt{3}}\sin\left(\frac{2\sqrt{2}}{\sqrt{3}}s\right), \frac{\sqrt{2}}{\sqrt{3}}\cos\left(\frac{1}{\sqrt{6}}s\right), \frac{\sqrt{2}}{\sqrt{3}}\sin\left(\frac{1}{\sqrt{6}}s\right)\right)$$

for all $s \in \mathbb{R}$. The curvatures are computed as $\kappa_\alpha = \frac{5}{\sqrt{2}}$ and $\tau_\alpha = \frac{2}{3}\varepsilon, \varepsilon = \pm 1$ [4]. Then from Corollary 2.2, $\alpha$ is not a Mannheim curve.